\documentstyle[11pt,twoside]{article}
\newtheorem{theorem}{\noindent\bf Theorem}[section]
\newtheorem{proposition}[theorem]{\noindent\bf Proposition}
\newtheorem{lemma}[theorem]{\noindent\bf Lemma}
\newtheorem{corollary}[theorem]{\noindent\bf Corollary}

\newtheorem{remark}[theorem]{\noindent\bf Remark}

\newtheorem{example}{\noindent\bf Example}

\newcommand{\e}{\hfill\blacksquare}
\input{amssym}
\begin{document}
\date{}
\title{{\Large\bf The multiplier algebra and
BSE-functions for\\ certain product of Banach algebras}}

\author{{\normalsize\sc Mehdi Nemati and Hossein Javanshiri}}

\maketitle
\normalsize

\begin{abstract}
In this paper, we characterize the (left) multiplier algebra of a
semidirect product algebra ${\mathcal A}={\mathcal B}\oplus
{\mathcal I}$, where ${\mathcal I}$ and ${\mathcal B}$ are closed
two-sided ideal  and closed subalgebra of ${\mathcal A}$,
respectively. As an application of this result  we investigate the
BSE-property of this class of Banach algebras. We then for two
commutative semisimple Banach algebras ${\mathcal A}$ and
${\mathcal B}$ characterize the BSE-functions on the carrier space
of ${\mathcal A}\times_\phi {\mathcal B}$, the  $\phi$-Lau product
of ${\mathcal A}$ and ${\mathcal B}$, in terms of those functions
on  carrier spaces of ${\mathcal A}$ and ${\mathcal B}$. We  also
prove that ${\mathcal A}\times_\phi {\mathcal B}$ is a BSE-algebra
if and only if both ${\mathcal A}$ and ${\mathcal B}$ are BSE.\\

{\bf Mathematics Subject Classification}: Primary 46J05, 46J10;
Secondary 46J99

{\bf Key words}: BSE-algebra, commutative Banach algebra,
   locally compact group, multiplier algebra

\end{abstract}


\section{Introduction}
Let $\mathcal A$  be a Banach algebra and suppose that ${\mathcal
I}$ and ${\mathcal B}$ are closed two-sided ideal  and closed
subalgebra of ${\mathcal A}$, respectively. Following Bade and
Dales \cite{bd, bd1} we say that $\mathcal A$ is a {\it semidirect
product} of $\mathcal B$ and $\mathcal I$ if ${\mathcal
A}={\mathcal B}\oplus {\mathcal I}$. Indeed, in this case the
product of two elements $(b, a)$ and $(b', a')$ of ${\mathcal
B}\oplus {\mathcal I}$ is given by
$$
(b, a)(b', a')=(bb', aa'+ba'+ab'),
$$
and this algebra endowed with the norm $\|(b, a)\|=\|b\|+\|a\|$ is
a Banach algebra. This notion was also studied by Berndt \cite{ob}
and Thomas \cite{thomas,thomas2} and they considered under what
conditions a commutative Banach algebra is the semidirect product
of a subalgebra and a principal ideal. We also note that the
algebra ${\mathcal B}\oplus{\mathcal I}$ is a  splitting extension
of $\mathcal B$ by $\mathcal I$. These products not only induce
new examples of Banach algebras which are interesting in its own
right but also they are known as a fertile source as examples or
counterexamples  in functional and abstract harmonic analysis; see
for example \cite{ll, cl}. An example of this type of algebras,
which is of special interest, is the {\it $\phi$-Lau product} of
two Banach algebras ${\mathcal A}$ and ${\mathcal B}$. In fact,
suppose that
 $\phi:{\mathcal B}\rightarrow {\mathcal A}$ is a  contractive algebra homomorphism. Then the { $\phi$-Lau product}
of  ${\mathcal A}$ and ${\mathcal B}$, denoted by ${\mathcal
A}\times_\phi {\mathcal B}$, is defined as the $\ell^1$-product
space ${\mathcal A}\times{\mathcal B}$ endowed with the norm
$\|(a,b)\|=\|a\|+\|b\|$ and the  product
$$
(a,b)(a',b')=(aa'+\phi(b)a'+a\phi(b'),bb')
$$
for all $(a, b), (a', b')\in {\mathcal A}\times{\mathcal B}$. It
is clear that with this norm and product, ${\mathcal
A}\times_{\phi}{\mathcal B}$ is a Banach algebra. Identifying
${\mathcal A}$ with ${\mathcal A}\times\{0\}$ and ${\mathcal B}$
with $\{0\}\times{\mathcal B}$, ${\mathcal A}$ is a closed ideal
of ${\mathcal A}\times_\phi {\mathcal B}$ and ${\mathcal B}$ is a
closed subalgebra. Therefore, ${\mathcal A}\times_{\phi}{\mathcal
B}$ is a semidirect product algebra of  ${\mathcal B}$ and
${\mathcal A}$. This type of product was first introduced by Bhatt
and Dabhi \cite{BD} for the case where ${\mathcal A}$ is
commutative and  was extended by the authors for the general case
\cite{se2, se4}.

In this paper, we characterize the (left) multiplier algebra of
${\mathcal B}\oplus {\mathcal I}$ in terms of the individual
(left) multiplier algebras of ${\mathcal B}$ and ${\mathcal I}$.
We then apply this result to investigate the BSE-property of this
class of Banach algebras. For two commutative semisimple Banach
algebras ${\mathcal A}$ and ${\mathcal B}$ we characterize the
BSE-functions on the carrier space $\Delta({\mathcal A}\times_\phi
{\mathcal B})$ of ${\mathcal A}\times_\phi {\mathcal B}$, in terms
of those functions on $\Delta({\mathcal A})$ and $\Delta({\mathcal
B})$. We also show that the $\phi$-Lau product of two commutative
semisimple Banach algebras ${\mathcal A}$ and ${\mathcal B}$ is
BSE if and only if both ${\mathcal A}$ and ${\mathcal B}$ are BSE.
 Finally, we use these results to study  the
BSE-property of certain Banach algebras related to a locally
compact group $G$.

\section{\large\bf (left) Multipliers on semidirect product algebras}

Let ${\mathcal C}$ be a Banach algebra and let  $X$ and $Y$ be two
Banach ${\mathcal C}$-bimodules. An operator $T: X\rightarrow Y$
is called a {\it left ${\mathcal C}$-module map} if $T(c\cdot x) =
c\cdot T(x)$ for all $c\in{\mathcal C}$ and $x\in X$. {\it Right
${\mathcal C}$-module} and {\it ${\mathcal C}$-bimodule maps}  are
defined similarly. We denote by $Hom_ {\mathcal C}(X,  Y)$ the
space of all bounded left ${\mathcal C}$-module maps from $ X$
into $Y$.  Define $LM({\mathcal C})$ to be the left  multiplier
algebra of ${\mathcal C}$. That is,
$$
LM({\mathcal C})=Hom_ {\mathcal C}({\mathcal C},  {\mathcal C}).
$$
Suppose now that $\mathcal A$  is the semidirect product algebra
${\mathcal B}\oplus {\mathcal I}$. Note that the dual ${\mathcal
A}^*$ of ${\mathcal A}$ can be identified with ${\mathcal
B}^{*}\times {\mathcal I}^{*}$ in the natural way given for each
$f\in {\mathcal I}^*$, $g\in {\mathcal B}^*$ and $(b,a)\in
{{\mathcal B}\oplus {\mathcal I}}$ by
$$\langle (g, f), (b,a) \rangle=g(b)+f(a).$$


We start the section with a key lemma.

\begin{lemma}\label{12}
Let $\mathcal A$  be the semidirect product algebra ${\mathcal
B}\oplus {\mathcal I}$.  Then $T\in LM({\mathcal A})$ if and only
if there exist some $R_{\mathcal I}\in Hom_{\mathcal B}({\mathcal
I}, {\mathcal I})$, $S_{\mathcal I}\in Hom_{\mathcal B}({\mathcal
B}, {\mathcal I})$, $T_{\mathcal B}\in LM({\mathcal B})$ and
$S_{\mathcal B}\in Hom_{\mathcal B}({\mathcal I}, {\mathcal B})$
such that for each $a, a'\in I$ and $b, b'\in {\mathcal B}$ we
have

\emph{(i)} $T((b, a)) = (S_{\mathcal B}(a) + T_{\mathcal B}(b),
R_{\mathcal I}(a) + S_{\mathcal I}(b))$.

\emph{(ii)} $R_{\mathcal I}(aa')=aR_{\mathcal I}(a')+aS_{\mathcal
B}(a')$.

\emph{(iii)} $R_{\mathcal I}(ab)=aS_{\mathcal I}(b)+aT_{\mathcal
B}(b)$.

\emph{(iv)} $S_{\mathcal B}(aa')=S_{\mathcal B}(ab)=0$.
\end{lemma}
{\noindent Proof.}
 Suppose that $T\in LM({\mathcal A})$. Then
 there exist bounded linear mappings $T_1: {\mathcal
A}\rightarrow {\mathcal B}$  and $T_2: {\mathcal A}\rightarrow
{\mathcal I}$  such that $T = (T_1, T_2)$. Let $S_{\mathcal B}(a)=
T_1((a, 0))$,
 $S_{\mathcal I}(b) =T_2((0, b))$, $T_{\mathcal B}(b) = T_1((0, b))$ and
 $R_{\mathcal I}(a) = T_2((a, 0))$ for all $a\in{\mathcal I}$ and $b\in {\mathcal B}$.
 Then trivially $S_{\mathcal B}$, $S_{\mathcal I}$, $T_{\mathcal B}$ and $R_{\mathcal I}$ are linear mappings satisfying (i). Moreover, for every
$a,a'\in {\mathcal I}$ and $b,b'\in {\mathcal B}$
\begin{eqnarray*}
(1)~~~~~~~~~~~~~~~~~~~~~~~\quad T{\big(}(b, a)( b',
a'){\big)}&=&T{\big(}(bb', aa'+b
a'+ab'){\big)}\\
&=&{\bigg(}S_{\mathcal B}(aa'+b a'+ab')+ T_{\mathcal B}(bb'),\\
&& R_{\mathcal I}(aa'+b a'+ab')+S_{\mathcal I}(bb'){\bigg)}
~~~~~~~~~~~~~~~~~~~~~~~~~~~~~~
\end{eqnarray*}
and
\begin{eqnarray*}
(2)~~~~~~~~~~~~~~~~~~~~\quad (b, a)T((b', a')) &=&(b,
a){\bigg(}S_{\mathcal B}(a')+T_{\mathcal B}(b'),
R_{\mathcal I}(a')+S_{\mathcal I}(b'){\bigg )} \\
&=&{\bigg(}b S_{\mathcal B}(a')+b T_{\mathcal B}(b'),
 aR_{\mathcal I}(a')+aS_{\mathcal I}(b')\\
&&+aS_{\mathcal B}(a')+aT_{\mathcal B}(b')+b R_{\mathcal I}(a')+b
S_{\mathcal I}(b') {\bigg)}~~~~~~~~~~~~~~~~~~~~~~~~~~~~~~~
\end{eqnarray*}
Thus $T\in LM({\mathcal A})$ if and only if $(1)$ and $(2)$
coincide; that is,
\begin{eqnarray*}
(3)~~~~~~~~~~~~~~~~~~~~~~~~~~~~~~~~\quad S_{\mathcal B}(aa'+b
a'+ab')+ T_{\mathcal B}(bb') &=&b S_{\mathcal B}(a')+b T_{\mathcal
B}(b'),~~~~~~~~~~~~~~~~~~~~~~~~~~
\end{eqnarray*}
and
\begin{eqnarray*}
(4)~~~~~~~~~~~~~~~~~~~\quad R_{\mathcal I}(aa'+b
a'+ab')+S_{\mathcal I}(bb') &=&aR_{\mathcal I}(a')+aS_{\mathcal
I}(b')+aS_{\mathcal B}(a')\\&&+aT_{\mathcal B}(b')+b R_{\mathcal
I}(a')+b S_{\mathcal I}(b').~~~~~~~~~~~~~~~~~~~~~~~~~~~~~~~~~~~~~~
\end{eqnarray*}
Therefore $T\in LM({\mathcal A})$ if and only if the equations
$(3)$ and $(4)$ are satisfied. Now a straightforward verification
shows that if $R_{\mathcal I}\in Hom_{\mathcal B}({\mathcal I},
{\mathcal I})$, $S_{\mathcal I}\in Hom_{\mathcal B}({\mathcal B},
{\mathcal I})$, $T_{\mathcal B}\in LM({\mathcal B})$ and
$S_{\mathcal B}\in Hom_{\mathcal B}({\mathcal I}, {\mathcal B})$
and the equalities (ii), (iii), (iv) are satisfied, then $(3)$ and
$(4)$ are valid. Applying $(3)$ and $(4)$ for suitable values of
$a,a', b, b'$ shows that $R_{\mathcal I}\in Hom_{\mathcal
B}({\mathcal I}, {\mathcal I})$, $S_{\mathcal I}\in Hom_{\mathcal
B}({\mathcal B}, {\mathcal I})$, $T_{\mathcal B}\in LM({\mathcal
B})$ and $S_{\mathcal B}\in Hom_{\mathcal B}({\mathcal I},
{\mathcal B})$ and the equalities (ii), (iii), (iv) are also
satisfied, as claimed.$\e$\\


For a right  Banach $\mathcal C$-module $X$ we denote by $\langle
X{\mathcal C}\rangle$ the closed linear span of the set
$X{\mathcal C}:=\{x\cdot c: c\in {\mathcal C}, x\in X\}$. As an
immediate consequence of Lemma \ref{12} we have the following
result.

\begin{corollary}\label{bi2}
Let $\mathcal A$  be the semidirect product algebra ${\mathcal
B}\oplus {\mathcal I}$ such that either $\langle{\mathcal
I}{\mathcal B}\rangle={\mathcal I}$ or $\langle{\mathcal
I}^2\rangle={\mathcal I}$.  Then $T\in LM({\mathcal A})$ if and
only if there exist $T_{\mathcal I}\in Hom_{\mathcal B}({\mathcal
I}, {\mathcal I})\cap LM({\mathcal I})$, $S_{\mathcal I}\in
Hom_{\mathcal B}({\mathcal B}, {\mathcal I})$ and $T_{\mathcal
B}\in LM({\mathcal B})$ such that $T_{\mathcal I}(ab)=aS_{\mathcal
I}(b)+aT_{\mathcal B}(b)$ and
$$
T((b, a)) = (T_{\mathcal B}(b), T_{\mathcal I}(a) + S_{\mathcal
I}(b))
$$
for all $a\in I$ and $b\in {\mathcal B}$.\\
\end{corollary}


\begin{remark}
{\rm Let $\mathcal A$  be a Banach algebra and suppose that
${\mathcal I}$ and ${\mathcal B}$ are closed two-sided ideal  and
closed subalgebra of ${\mathcal A}$, respectivly. Then for each
$\varphi\in \Delta({\mathcal I})$ there is a unique
$\psi_\varphi\in \Delta({\mathcal B})\cup\{0\}$ such that
$\varphi(ab)=\varphi(ba)=\varphi(a)\psi_\varphi(b)$ for all
$a\in{\mathcal I}$ and $b\in{\mathcal B}$. Indeed,
$\psi_\varphi\in \Delta({\mathcal B})\cup\{0\}$ is defined by
$\psi_\varphi(b):=\varphi(ba_0)$  for all $b\in{\mathcal B}$,
where $a_0\in {\mathcal I}$ is any element with $\varphi(a_0)=1$;
to see this, note that for each $a\in{\mathcal I}$ and
$b\in{\mathcal B}$,
\begin{eqnarray*}
\varphi(ab)=\varphi(aba_0)=\varphi(a)\varphi(ba_0)=\varphi(a)\psi_\varphi(b).
\end{eqnarray*}
That  $\psi_\varphi$  is unique follows trivially. We note that if
$\langle{\mathcal I}{\mathcal B}\rangle={\mathcal I}$, then
$\psi_\varphi\neq 0$ for all $\varphi\in \Delta(I)$.}
\end{remark}


\begin{proposition}\label{ib}
Let $\mathcal A$  be the semidirect product algebra ${\mathcal
B}\oplus {\mathcal I}$ and let
\begin{eqnarray*}
E:=\{(\psi_\varphi, \varphi): \varphi\in \Delta({\mathcal
I})\}\quad{\rm and}\quad F:=\{(\psi, 0): \psi\in \Delta({\mathcal
B})\}.
\end{eqnarray*}
Then $E$ and $F$ are disjoint and  $\Delta({\mathcal A})=E\cup F$.
\end{proposition}
{\noindent Proof.} It is clear that $E\cap F=\emptyset$ and $E\cup
F\subseteq \Delta({\mathcal A})$. For the converse, suppose that
$(\psi, \varphi)\in \Delta({\mathcal A})$. Then, for each $(b, a),
(b', a')\in {\mathcal A}$ we have
$$
\langle (\psi, \varphi), (b, a)(b', a')\rangle =\langle (\psi,
\varphi), (b, a)\rangle \langle (\psi, \varphi), (b', a')\rangle.
$$
This implies that
$$\psi(bb')+\varphi(ab'+ba'+aa')=\psi(b)\psi(b')+\psi(b)\varphi(a')+\varphi(a)\psi(b')+\varphi(a)\varphi(a').$$
Taking $b=b'=0$, it follows that
$\varphi(aa')=\varphi(a)\varphi(a')$ and therefore $\varphi\in
\Delta({\mathcal I})\cup\{0\}$. Similarly, we can see that
$\psi\in \Delta({\mathcal B})\cup\{0\}$. Now, if $a=b'=0$, then we
have $\varphi(ba')=\psi(b)\varphi(a')$, similarly
$\varphi(ab')=\varphi(a)\psi(b')$ for all $a,a'\in {\mathcal I}$
and $b, b'\in {\mathcal B}$.  The equality $\varphi=0$ implies
that $\psi\neq 0$. If
 $\varphi\neq 0$, then $\psi=\psi_\varphi$ by the above remark.$\e$\\


\section{\large\bf The BSE-property of semidirect product algebras}

A commutative Banach algebra ${\mathcal C}$ is called without
order if for each $c\in{\mathcal C}$, $c{\mathcal C}=\{0\}$
implies $c=0$. For example, if ${\mathcal C}$ is a commutative
semisimple Banach algebra, then it is without order. Let
${\mathcal C}$ be a commutative Banach algebra with carrier space
$\Delta({\mathcal C})$ and let ${\mathcal C}^*$ denote the dual
space of ${\mathcal C}$. A continuous complex-valued function
$\sigma$ on $\Delta({\mathcal C})$ is said to satisfy {\it the
Bochner-Schoenberg-Eberlein {\rm(}BSE{\rm)} inequality} if there
exists a constant $C>0$ such that for any
  $\varphi_1,..., \varphi_n\in\Delta({\mathcal C})$ and  $c_1,..., c_n\in {\Bbb{C} }$ the
inequality
$$
\left|\sum_{j=1}^n c_j\sigma(\varphi_j)\right|\leq C\cdot
\left\|\sum_{j=1}^n c_j\varphi_j\right\|_{{\mathcal C}^*}
$$
holds. Let $C_{BSE}(\Delta({\mathcal C}))$ denote the set of  all
continuous complex-valued functions on $\Delta({\mathcal C})$
satisfying the BSE-inequality. The BSE-norm of $\sigma$ denoted by
$\|\sigma\|_{BSE}$, is defined to be the infimum of all such $C$.
Takahasi and Hatori \cite{hat1} showed that under this norm
$C_{BSE}(\Delta({\mathcal C}))$ is a commutative semisimple Banach
algebra. A linear operator $T$ on ${\mathcal C}$ is called a
multiplier if it satisfies $cT(b)=T(c)b$ for all $b, c\in
{\mathcal C}$. Suppose that $M({\mathcal C})$ denotes the space of
all  multiplier of the commutative Banach algebra ${\mathcal C}$
which is a unital commutative Banach algebra. Recall that  for
each $T\in M({\mathcal C})$ there exists a unique continuous
function $\widehat{T}$ on $\Delta({\mathcal C})$ such that
$\widehat{T(c)}(\varphi)=\widehat{T}(\varphi)\widehat{c}(\varphi)$
for all $c\in{\mathcal C}$ and $\varphi\in\Delta({\mathcal C})$;
see  \cite[Theorem 1.2.2.]{lar}. A commutative Banach algebra
${\mathcal C}$ without order  is called a {\it BSE-algebra} (or is
said to have the {\it BSE-property}) if
$$
C_{BSE}(\Delta({\mathcal C}))=\widehat{M({\mathcal C})}.
$$
The concept of BSE-property  has been first introduced and studied
by Takahasi and Hatori \cite{hat1} and later  by several authors
in various classes of commutative Banach algebras; see also
\cite{ino2, ino1, iz,  lash2, kani3, ulg, hat2}.

 A bounded net $(e_\alpha)_\alpha$ in a Banach
algebra  ${\mathcal C}$ is called a {\it $\Delta$-weak bounded
approximate identity} if it satisfies
$\varphi(e_\alpha)\rightarrow 1$  for all $\varphi\in
\Delta({\mathcal C})$. Such approximate identities were studied in
\cite{jon}. It was shown in \cite[Corollary 5]{hat1} that a
commutative Banach algebra ${\mathcal C}$ has a $\Delta$-weak
bounded approximate identity if and only if $\widehat{M({\mathcal
C})} \subseteq C_{BSE}(\Delta({\mathcal C}))$.

 In the sequel, we assume that ${\mathcal A}$ is a commutative semidirect product algebra of semisimple closed subalgebra ${\mathcal B}$ and  semisimple closed ideal ${\mathcal I}$.


\begin{corollary}\label{prod}
Let $\mathcal A$  be the commutative semidirect product algebra
${\mathcal B}\oplus {\mathcal I}$.  Then $T\in M({\mathcal A})$ if
and only if there exist $T_{\mathcal I}\in Hom_{\mathcal
B}({\mathcal I}, {\mathcal I})\cap M({\mathcal I})$, $S_{\mathcal
I}\in Hom_{\mathcal B}({\mathcal B}, {\mathcal I})$ and
$T_{\mathcal B}\in M({\mathcal B})$ such that $T_{\mathcal
I}(ab)=aS_{\mathcal I}(b)+aT_{\mathcal B}(b)$ and
$$
T((b, a)) = (T_{\mathcal B}(b), T_{\mathcal I}(a) + S_{\mathcal
I}(b))
$$
for all $a\in I$ and $b\in {\mathcal B}$.
\end{corollary}
{\noindent Proof.} Suppose that $T\in M({\mathcal A})$.  By Lemma
\ref{12}, there exist some $R_{\mathcal I}\in Hom_{\mathcal
B}({\mathcal I}, {\mathcal I})$, $S_{\mathcal I}\in Hom_{\mathcal
B}({\mathcal B}, {\mathcal I})$, $T_{\mathcal B}\in M({\mathcal
B})$ and $S_{\mathcal B}\in Hom_{\mathcal B}({\mathcal I},
{\mathcal B})$ such that for each $a\in I$ and $b\in {\mathcal B}$
we have
$$T((b, a)) = (S_{\mathcal B}(a) + T_{\mathcal B}(b),
R_{\mathcal I}(a) + S_{\mathcal I}(b)).$$ Moreover, $S_{\mathcal
B}(ba)=0$. Since $S_{\mathcal B}\in Hom_{\mathcal B}({\mathcal I},
{\mathcal B})$ and ${\mathcal B}$ is without order, it follows
that $S_{\mathcal B}=0$ and $R_{\mathcal I}\in M({\mathcal I})$.
The converse is clear.$\e$\\


\begin{proposition}\label{sub}
Let $\mathcal A$  be the semidirect product algebra ${\mathcal
B}\oplus {\mathcal I}$ such that $\langle{\mathcal I}{\mathcal
B}\rangle={\mathcal I}$.  Then ${\mathcal B}$ is a BSE-algebra if
${\mathcal A}$ is so.
\end{proposition}
{\noindent Proof.} Suppose that $\rho\in C_{BSE}(\Delta ({\mathcal
B}))$. Define a function $\sigma$ on $\Delta({\mathcal A})$ by
$$
\sigma(\psi, 0)=\rho(\psi),\quad \sigma(\psi_\varphi,
\varphi)=\rho(\psi_\varphi)
$$
for all $\psi\in \Delta({\mathcal B})$ and $\varphi\in
\Delta({\mathcal I})$. Then  $\sigma$ is continuous since
$\langle{\mathcal I}{\mathcal B}\rangle={\mathcal I}$ and so
$\psi_\varphi\neq 0$ for all $\varphi\in \Delta(I)$. Moreover,
since $\rho\in C_{BSE}(\Delta ({\mathcal B}))$,  by \cite[Theorem
4(i)]{hat1} there is a bounded net $(b_\alpha)_\alpha$ in
${\mathcal B}$ such that $\widehat{b_\alpha}(\psi)\rightarrow
\sigma(\psi)$ for all $\psi\in \Delta({\mathcal B})$. Now, if we
consider $(b_\alpha, 0)_\alpha$ as a net in ${\mathcal A}$, then
$$
\widehat{(b_\alpha, 0)}(\psi, 0)\rightarrow
\sigma(\psi)=\sigma(\psi, 0)
$$
for all $\psi\in \Delta({\mathcal B})$ and
$$
\widehat{(b_\alpha, 0)}(\psi_\varphi, \varphi)\rightarrow
\rho(\psi_\varphi)=\sigma(\psi_\varphi, \varphi)
$$
for all $\varphi\in \Delta({\mathcal I})$. Thus, $\sigma\in
C_{BSE}(\Delta ({\mathcal A}))$ and therefore, there is some $T\in
M({\mathcal A})$ such that $\widehat{T}=\sigma$. By Corollary
\ref{prod}, there exist some $T_{\mathcal I}\in Hom_{\mathcal
B}({\mathcal I}, {\mathcal I})$, $S_{\mathcal I}\in Hom_{\mathcal
B}({\mathcal B}, {\mathcal I})$, $T_{\mathcal B}\in M({\mathcal
B})$ such that for each $a\in I$ and $b\in {\mathcal B}$ we have
$$
T((b, a)) = (T_{\mathcal B}(b), T_{\mathcal I}(a) + S_{\mathcal
I}(b))
$$
On the other hand, $\widehat{T(b, 0)}(\psi, 0)=\widehat{T}(\psi,
0)\widehat{(b, 0)}(\psi, 0)$ for all $\psi\in \Delta({\mathcal
B})$ and $b\in {\mathcal B}$. Thus,
$$
\psi(T_{\mathcal B}(b))=\widehat{T}(\psi, 0)\psi(b).
$$
It follows that $\widehat{T_{\mathcal
B}(b)}(\psi)=\widehat{T}(\psi, 0)\widehat{b}(\psi)$. Moreover,
$\widehat{T_{\mathcal B}(b)}(\psi)=\widehat{T_{\mathcal
B}}(\psi)\widehat{b}(\psi)$ for all $b\in {\mathcal B}$.
Therefore, $\widehat{T}((\psi, 0))=\widehat{T_{\mathcal B}}(\psi)$
and consequently,
$$
\rho(\psi)=\sigma(\psi, 0)=\widehat{T}(\psi,
0)=\widehat{T_{\mathcal B}}(\psi)
$$
for all $\psi\in \Delta({\mathcal B})$. This means that $\rho\in
\widehat{M(\mathcal B)}$; that is, $C_{BSE}(\Delta({\mathcal
B}))\subseteq\widehat{M(\mathcal B)}$. The reverse inclusion
follows from this fact that if ${\mathcal A}$ has a $\Delta$-weak
bounded approximate identity, then so does ${\mathcal
B}$.$\e$\\


Now, for two Banach algebras ${\mathcal A}$ and ${\mathcal B}$ we
turn our attention to the BSE-property of direct sum algebra
${\mathcal A}\times_0 {\mathcal B}$. Recall that ${\mathcal
A}\times_0 {\mathcal B}$ is equipped with the usual direct product
multiplication and the norm defined by $\|(a, b)\|=\|a\|+\|b\|$
for all $a\in {\mathcal A}$ and $b\in{\mathcal B}$. Also, it is
clear that $\Delta({\mathcal A}\times_0 {\mathcal B})=E\cup F$,
where
$$
E=\{(\varphi, 0): \varphi\in \Delta({\mathcal A})\}\quad
\hbox{and}\quad F=\{(0, \psi): \psi\in \Delta({\mathcal B})\}.
$$
We have the following corollary which was obtained independently
in \cite{lash2} through different methods. For sake of
completeness and as an application of Corollary \ref{prod}, we
will give a  proof of the corollary.


\begin{corollary}\label{tim2}
Let $\mathcal A$ and $\mathcal B$  be  two commutative semisimple
Banach algebras. Then  ${\mathcal A}\times_0 {\mathcal B}$ is a
BSE-algebra if and only if ${\mathcal A}$ and  ${\mathcal B}$ are
BSE-algebras.
\end{corollary}
{\noindent Proof.} It is easily verified that
$C_{b}(\Delta({\mathcal A}\times_0 {\mathcal
B}))=C_{b}(\Delta({\mathcal A}))\times_0 C_b(\Delta({\mathcal
B}))$. Indeed,  $\sigma\in C_{b}(\Delta({\mathcal A}\times_0
{\mathcal B}))$ if and only if there are   $\tau\in
C_{b}(\Delta({\mathcal A}))$ and $\rho\in C_{b}(\Delta({\mathcal
B}))$ such that $\sigma=(\tau, \rho)$; that is,
$$
\sigma(\varphi, 0)=\tau(\varphi)\quad{\rm and}\quad\sigma(0,
\psi)=\rho(\psi)
$$
for all $\varphi\in \Delta({\mathcal A})$ and $\psi\in
\Delta({\mathcal B})$. Moreover,
$$C_{BSE}(\Delta({\mathcal A}\times_0 {\mathcal B}))=
C_{b}(\Delta({\mathcal A}\times_0 {\mathcal B})) \cap({\mathcal
A}^{**}\times_0 {\mathcal B}^{**})|_{\Delta({\mathcal A}\times_0
{\mathcal B})},$$ by  \cite[Theorem 4]{hat1}. This shows that
$C_{BSE}(\Delta({\mathcal A}\times_0 {\mathcal
B}))=C_{BSE}(\Delta({\mathcal A}))\times_0
C_{BSE}(\Delta({\mathcal B}))$. On the other hand, by Corollary
\ref{prod}  we have $M({\mathcal A}\times_0 {\mathcal
B})=M({\mathcal A})\times_0 M(\mathcal B)$. It is also easy to
check that $$\widehat{M({\mathcal A}\times_0 {\mathcal
B})}=\widehat{M({\mathcal A})}\times_0\widehat{ M(\mathcal B)}.$$
The rest of the proof is straightforward.$\e$\\




\section{\large\bf\bf The BSE-property of $\phi$-Lau product algebras}

Let ${\mathcal A}$, ${\mathcal B}$ and $\phi$ be as in the
introduction. If $\phi=0$, then the algebra ${\mathcal
A}\times_\phi {\mathcal B}$ coincides with the direct sum algebra
${\mathcal A}\times_0{\mathcal B}$. We also recall from
\cite[Theorem 2.2]{se2} that $\Delta({\mathcal
A}\times_{\phi}{\mathcal B})=E\cup F$, where
\begin{eqnarray*}
E=\{(\varphi,\varphi\circ\phi):~~\varphi\in\Delta({\mathcal
A})\}\quad {\rm and}\quad F=\{(0,\psi):~~\psi\in\Delta({\mathcal
B})\}.
\end{eqnarray*}
 In the sequel for two  commutative semisimple Banach algebras ${\mathcal A}$ ad ${\mathcal B}$ we characterize the BSE-functions on $\Delta({\mathcal
A}\times_\phi {\mathcal B})$ in terms of those functions on
$\Delta({\mathcal A})$ and $\Delta({\mathcal B})$.


\begin{lemma}\label{th0}
 Let $\mathcal A$ and $\mathcal B$  be commutative semisimple Banach algebras and let $\phi:{\mathcal B}\rightarrow {\mathcal A}$ be
 a contractive algebra homomorphism with dense range. Then the following statements hold.

 \emph{(i)} Let $\sigma\in C_{BSE}(\Delta ({\mathcal A}\times_\theta {\mathcal B}))$ and define functions $\tau$ on $\Delta ({\mathcal A})$  and $\rho$ on $\Delta ({\mathcal B})$ by
 $$
 \tau(\varphi)=\sigma(\varphi, \varphi\circ\phi)-\sigma(0, \varphi\circ\phi) \quad {\rm and}\quad \rho(\psi)=\sigma(0, \psi)
 $$
for all $\varphi\in \Delta({\mathcal A})$ and $\psi\in
\Delta({\mathcal B})$.  Then $\tau\in C_{BSE}(\Delta ({\mathcal
A}))$ and  $\rho\in C_{BSE}(\Delta ({\mathcal B}))$.

 \emph{(ii)} Let $\tau\in C_{BSE}(\Delta ({\mathcal A}))$ and $\rho\in C_{BSE}(\Delta ({\mathcal B}))$ and  define a function $\sigma$ on $\Delta ({\mathcal A}\times_\theta {\mathcal B})$ by
 $$\sigma(\varphi, \varphi\circ\phi)=\tau(\varphi)+\rho(\varphi\circ\phi) \quad {\rm and}\quad \sigma(0, \psi)=\rho(\psi)$$ for all $\varphi\in \Delta({\mathcal A})$ and $\psi\in \Delta({\mathcal B})$. Then $\sigma\in C_{BSE}(\Delta ({\mathcal A}\times_\phi {\mathcal B}))$.
\end{lemma}
{\noindent Proof.} (i). First note that the density of
$\phi({\mathcal B})$ in ${\mathcal A}$ implies that
$\varphi\circ\phi\neq 0$ for all $\varphi\in \Delta({\mathcal
A})$, and therefore $\tau$ is well defined. It is also trivial
that $\tau$ and $\rho$ are continuous. Moreover,  by \cite[Theorem
4(i)]{hat1} there is a bounded net $(a_\alpha, b_\alpha)_\alpha$
in ${\mathcal A}\times_\phi {\mathcal B}$ such that $\|(a_\alpha,
b_\alpha)\|\leq\|\sigma\|_{BSE}$ and
 $\widehat{(a_\alpha, b_\alpha)}(\gamma)\rightarrow \sigma(\gamma)$ for all $\gamma\in \Delta({{\mathcal A}\times_\phi {\mathcal B}})$. Therefore, for each $\psi\in\Delta({\mathcal B})$,
 $$\widehat{b_\alpha}(\psi)=\widehat{(a_\alpha, b_\alpha)}(0, \psi)\rightarrow \sigma(0, \psi) =\rho(\psi)$$ and for each $\varphi\in\Delta({\mathcal A})$ we have
 $$
 \widehat{a_\alpha}(\varphi)=\widehat{(a_\alpha, b_\alpha)}(\varphi, \varphi\circ\phi)-\widehat{b_\alpha}(\varphi\circ\phi)\rightarrow \sigma(\varphi, \varphi\circ\phi)-\sigma(0, \varphi\circ\phi)=\tau(\varphi).
 $$
 These show    that $\tau\in C_{BSE}(\Delta ({\mathcal A}))$ and
  $\rho\in C_{BSE}(\Delta ({\mathcal B}))$. Finally,  since
  $\|b_\alpha\|+ \|a_\alpha\|= \|(a_\alpha, b_\alpha)\|\leq \|\sigma\|_{BSE}$, it follows from \cite[Remark on p. 154]{hat1} that $\|\tau\|_{BSE}+\|\rho\|_{BSE}\leq\|\sigma\|_{BSE}$.

(ii). As it was mentioned  above
$\varphi\circ\phi\in{\Delta({\mathcal B})}$ for all $\varphi\in
\Delta({\mathcal A})$ and consequently $\sigma$ is well defined.
Furthermore, Since $E$ and $F$ are closed in $\Delta ({\mathcal
A}\times_\phi {\mathcal B})$, it follows that $\sigma$ is
continuous.   By assumption and \cite[Theorem 4(i)]{hat1} there
exist  bounded nets $(a_\alpha)_\alpha$ and $(b_\beta)_\beta$ in
${\mathcal A}$ and ${\mathcal B}$, respectively with $\|
a_\alpha\|\leq \|\tau\|_{BSE}$ and $\|b_\beta\|\leq
\|\rho\|_{BSE}$ such that
 $\widehat{a_\alpha}(\varphi)\rightarrow \tau(\varphi)$ and $\widehat{b_\beta}(\psi)\rightarrow \rho(\psi)$ for all $\varphi\in \Delta ({\mathcal A})$ and $\psi\in \Delta ({\mathcal B})$. If we consider the  bounded net $(a_\alpha, b_\beta)_{\alpha, \beta}$ in ${\mathcal A}\times_\phi {\mathcal B}$, then
 $$
 \widehat{(a_\alpha, b_\beta)}(\varphi, \varphi\circ\phi)= \widehat{a_\alpha}(\varphi)+ \widehat{b_\beta}(\varphi\circ \phi)\rightarrow \tau(\varphi)+\rho(\varphi\circ\phi)=\sigma(\varphi, \varphi\circ\phi)
 $$
 for all $(\varphi, \varphi\circ\phi)\in E$ and
$$
 \widehat{(a_\alpha, b_\beta)}(0, \psi)= \widehat{b_\beta}(\psi)\rightarrow \rho(\psi)=\sigma(0, \psi)
 $$
for all $(0, \psi)\in F$. Therefore, $\widehat{(a_\alpha,
b_\beta)}(\gamma)\rightarrow \sigma(\gamma)$ for all
$\gamma\in\Delta ({\mathcal A}\times_\phi {\mathcal B})$. Hence,
$\sigma\in C_{BSE}(\Delta ({\mathcal A}\times_\phi{\mathcal B}))$
and
 $\|\sigma\|_{BSE}\leq \|\tau\|_{BSE}+\|\rho\|_{BSE}$.$\e$\\


Now, assume that $\mathcal A$ and $\mathcal B$  are commutative
semisimple
 Banach algebras and let $\phi:{\mathcal B}\rightarrow
{\mathcal A}$ be a contractive algebra homomorphism with dense
range.  Suppose that $\Gamma=\phi^*|_{\Delta({\mathcal A})}:
\Delta({\mathcal A})\rightarrow \Delta({\mathcal B})$. Then the
map $\widetilde{\phi}:C_{BSE}({\Delta({\mathcal B})})\rightarrow
C_{BSE}({\Delta({\mathcal A})})$ defined by
$\widetilde{\phi}(\sigma)=\sigma\circ \Gamma$ is a contractive
algebra homomorphism.

\begin{theorem}
Let $\phi:{\mathcal B}\rightarrow {\mathcal A}$ be a
 contractive algebra homomorphism with dense range. Then the map $\Theta: (\tau, \rho)\rightarrow \sigma$ of Lemma \ref{th0} is an isometric isomorphism from
 $$
 C_{BSE}({\Delta({\mathcal B})})\times_{\widetilde{\phi}} C_{BSE}({\Delta({\mathcal A})})
 $$
 onto $C_{BSE}({\Delta({\mathcal A}\times_\phi{\mathcal B})})$.
\end{theorem}
{\noindent Proof.} It follows from Lemma \ref{th0} that $\Theta$
is bijective and isometric. It suffices to observe that $\Theta$
is also an algebra isomorphism. Indeed, given $\tau_i\in
C_{BSE}({\Delta({\mathcal A})})$ and $\rho_i\in
C_{BSE}({\Delta({\mathcal B})})$, $i=1, 2$. Then  for each
$\varphi\in \Delta({\mathcal A})$,
\begin{eqnarray*}
\langle \Theta((\tau_1, \rho_1)(\tau_2, \rho_2)), (\varphi, \varphi\circ\phi)\rangle&=& \langle \Theta(\tau_1\tau_2+\widetilde{\phi}(\rho_1)\tau_2+\tau_1\widetilde{\phi}(\rho_2), \rho_1\rho_2), (\varphi, \varphi\circ\phi)\rangle\\
&=& [\tau_1(\varphi)+\rho_1(\varphi\circ\phi)][\tau_2(\varphi)+\rho_2(\varphi\circ\phi)]\\
&=& \langle\Theta(\tau_1, \rho_1)\Theta(\tau_2, \rho_2), (\varphi,
\varphi\circ\phi)\rangle.
\end{eqnarray*}
Similarly we can show that $\langle \Theta((\tau_1,
\rho_1)(\tau_2, \rho_2)), (0, \psi)\rangle=\langle\Theta(\tau_1,
\rho_1)\Theta(\tau_2, \rho_2), (0, \psi)\rangle$  for all $\psi\in
\Delta({\mathcal B})$ which completes the proof.$\e$\\


Let $\mathcal A$ and $\mathcal B$  be commutative semisimple
Banach algebras and let $\phi:{\mathcal B}\rightarrow {\mathcal
A}$ be a contractive algebra homomorphism. Then it is easy to
check that the map $\Phi: {\mathcal A}\times_0 {\mathcal
B}\rightarrow {\mathcal A}\times_\phi {\mathcal B}$ defined by
$$
\Phi(a, b)=(a-\phi(b), b)\quad ((a, b)\in {\mathcal A}\times_0
{\mathcal B} )
$$
is an algebra isomorphism and
$$
\|(a-\phi(b),
b)\|=\|a-\phi(b)\|+\|b\|\leq\|a\|+(\|\phi\|+1)\|b\|\leq
(\|\phi\|+1)\|(a, b)\|
$$
Moreover, Since both algebras are semisimple, it follows that
$\Phi$ is also a topological isomorphism. Therefore, there is a
homeomorphism $\widehat{\Phi}$ between $\Delta({\mathcal
A}\times_0 {\mathcal B})$ and $\Delta({\mathcal
A}\times_\phi{\mathcal B})$ given by
$$
\widehat{\Phi}(\varphi, \varphi\circ \phi)=(\varphi,
0),\quad\widehat{\Phi}(0, \psi)=(0, \psi)
$$
for all $\varphi\in \Delta({\mathcal A})$ and $\psi\in
\Delta({\mathcal B})$. As an application of Corollary \ref{tim2}
we have the following result characterizing the BSE-property of
${\mathcal
 A}\times_\phi{\mathcal B}$.


\begin{theorem}
Let $\mathcal A$ and $\mathcal B$  be  two commutative semisimple
Banach algebras. Then  ${\mathcal A}\times_\phi {\mathcal B}$ is a
BSE-algebra if and only if ${\mathcal A}$ and  ${\mathcal B}$ are
BSE-algebra.
\end{theorem}
{\noindent Proof.} By Corollary \ref{tim2} it suffices to show
that ${\mathcal A}\times_\phi  {\mathcal B}$ is a BSE-algebra if
and only if ${\mathcal A}\times_0 {\mathcal B}$ is a BSE-algebra.
In fact, suppose that ${\mathcal A}\times_\phi  {\mathcal B}$ is a
BSE-algebra and let
 ${\sigma}\in  C_{BSE}(\Delta({\mathcal A}\times_0  {\mathcal B}))$. Define $\tau$ on $\Delta({\mathcal A}\times_\phi {\mathcal B})$ by
 $$
 \tau(\varphi, \varphi\circ\phi)=\sigma(\varphi, 0),\quad\tau(0, \psi)=\sigma(0, \psi)
 $$
 for all $\varphi\in \Delta({\mathcal A})$ and $\psi\in \Delta({\mathcal B})$. Since $\tau=\sigma\circ {\Phi}^*$, it follows that $\tau\in  C_{BSE}(\Delta({\mathcal A}\times_\phi  {\mathcal B}))$, where $\Phi$ is the above isomorphism between ${\mathcal A}\times_0 {\mathcal B}$ and ${\mathcal A}\times_\phi  {\mathcal B}$.
 By assumption  there is some $T\in M({\mathcal A}\times_\phi  {\mathcal B})$ such that $\widehat{T}=\tau$. Now,  consider the map $S=\Phi^{-1} T\Phi: {\mathcal A}\times_0  {\mathcal B}\rightarrow {\mathcal A}\times_0  {\mathcal B}$. Then it is trivial that $S\in M({\mathcal A}\times_0  {\mathcal B})$. We show that $\widehat{S}=\sigma$. Indeed, for each $\varphi\in \Delta({\mathcal A})$, $a\in {\mathcal A}$ and $b\in{\mathcal B}$ we have
 \begin{eqnarray*}
 \widehat{S}(\varphi, 0)\widehat{(a, b)}(\varphi, 0)&=&\widehat{S(a, b)}(\varphi, 0)\\
 &=&\langle \Phi^{-1}T\Phi(a, b), (\varphi, 0)\rangle\\
 &=&\langle T(a-\phi(b), b), {\Phi^{-1}}^*(\varphi, 0)\rangle\\
 &=&\langle T(a-\phi(b), b), (\varphi, \varphi\circ\phi)\rangle\\
 &=&\widehat{\tau}(\varphi, \varphi\circ\phi)\widehat{a}(\varphi)\\
 &=&{\sigma}(\varphi, 0)\widehat{(a, b)}(\varphi, 0).
 \end{eqnarray*}
Therefore, $\widehat{S}(\varphi, 0)=\sigma(\varphi, 0)$.
Similarly, we can show that $\widehat{S}(0, \psi)=\sigma(0, \psi)$
for all $\psi\in \Delta({\mathcal B})$. This shows that
$C_{BSE}(\Delta({\mathcal A}\times_0  {\mathcal
B}))\subseteq\widehat{M({\mathcal A}\times_0 {\mathcal B})}$. For
the reverse inclusion  note that if $(u_\alpha, v_\alpha)$ is a
$\Delta$-weak bounded approximate identity for ${\mathcal
A}\times_\phi  {\mathcal B}$, then $\Phi^{-1}((u_\alpha,
v_\alpha))$ is  a  $\Delta$-weak bounded approximate identity for
${\mathcal A}\times_0  {\mathcal B}$, which implies that
${\mathcal A}\times_0 {\mathcal B}$ is a BSE-algebra. Similarly,
we can show that the converse is also true.  $\e$\\



Recall that a locally compact group $G$ is called {\it amenable}
if there exists a translation-invariant mean on $L^\infty(G)$.
Now, let  $A(G)$ be the Fourier algebra  of $G$, introduced by
Eymard \cite{eymard}. Then $A(G)$ is always an ideal in the
Fourier-Stieltjes algebra $B(G)$ and $M(A(G))=B(G)$ when $G$ is
amenable. Furthermore, the carrier space of $A(G)$ as a
commutative Banach algebra is homeomorphic to the topological
space $G$ where each $x\in G$ is mapped to $\varphi_x$, the point
evaluation at $x$. Thus, in the following we identify
$\Delta(A(G))$ with $G$.

The following result was obtained independently in a different way
in \cite{ulg}.

\begin{theorem}
Let $G$ be an amenable locally compact group and let $N$ be a
closed normal subgroup of G such that $G/N$ is compact. Suppose
that  $q : G\rightarrow G/N$ is the quotient map, and let
$$B(G) = [B(G/N) \circ q] \times_0 A(G).$$
Then $B(G)$ is a BSE-algebra.
\end{theorem}
{\noindent Proof.} Since $G/N$  is compact, it follows from
\cite[Corollary 2.26(3)]{eymard} that the map $u\rightarrow u\circ
q$ is an isometric Banach algebra isomorphism between
$A(G/N)=B(G/N)$ and its image $A(G/N))\circ q$. Since both $G$ and
$G/N$ are amenable groups, it follows from Corollary \ref{tim2}
and \cite[Theorem 5.1]{ulg} that $B(G)$ is a
BSE-algebra.$\e$\\


The following  example shows that $B(G)$ can be a BSE-algebra when
$G$ is noncompact. Therefore, we do not have the dual version of
the result on BSE-property of $M(G)$.

\begin{example}
{\rm Let $n\geq 1$ be an integer, $p$ a prime number, ${\Bbb Q}_p$
the field of $p$-adic numbers, and ${\Bbb Z}_p$ the ring of
$p$-adic integers. Then the semidirect product $G_{p,n} = {\Bbb
Q}^n_p \ltimes GL(n,{\Bbb Z}_p)$  is a noncompact amenable locally
compact group. Runde and Spronk \cite{rs} showed that
$$B(G_{p,n}) = [A(GL(n,{\Bbb Z}_p) )\circ q] \times_0 A(G_{p,n}).$$
Since $GL(n,{\Bbb Z}_p)$ and $G_{p,n}$ are compact and amenable
groups, respectively, it follows from above theorem that
$B(G_{p,n})$ is a BSE-algebra.}
\end{example}

\vspace{3mm}

\noindent {\sc Mehdi Nemati}\\
Department of Mathematical Sciences, Isfahan Uinversity of
Technology, Isfahan 84156-83111, Iran;\\
 and\\
School of Mathematics, Institute for Research in Fundamental
Sciences (IPM), P.O. Box: 19395–5746, Tehran, Iran\\
E-mail: m.nemati@cc.iut.ac.ir\\

\noindent {\sc Hossein Javanshiri}\\
Department of Mathematics, Yazd University,
P.O. Box: 89195-741, Yazd, Iran\\
{ E-mail}: h.javanshiri@yazd.ac.ir\\

\end{document}